\documentclass[12pt]{article}
\usepackage{amssymb,latexsym}
\usepackage{amsfonts}
\title{Generalized Umemura polynomials and Hirota--Miwa equations}
\author{Anatol N. Kirillov and Makoto Taneda}
\date{}

\begin{document}

\newcommand{\bbb}{{\bf b}}%
\newcommand{\C}{{\mathbb C}}%
\newcommand{\HH}{{\cal H}}%
\newcommand{\dd}{{\bf d}}
\newcommand{\cc}{{\bf c}}
\newcommand{\odd}{{\rm odd}}
\newcommand{\even}{{\rm even}}
\newcommand{\Split}{{\bf Split}}

\newcounter{newsection}
\renewcommand{\thesection}{\S \arabic{section}.}
\renewcommand{\thesubsection}{\thesection.\arabic{subsection}}
\renewcommand{\theequation}{\thenewsection.\arabic{equation}}
\renewcommand
{\thenewsection}{\setcounter{newsection}{\value{section}}\arabic{newsection}}
\newtheorem{Proposition}{Proposition}
\newtheorem{Definition}{Definition}
\newtheorem{Theorem}{Theorem}
\newtheorem{Remarks}{Remarks}
\newtheorem{Corollary}{Corollary}
\newtheorem{Conjecture}{Conjecture}
\newtheorem{Lemma}{Lemma}
\maketitle
\begin{abstract}
We introduce and study generalized Umemura polynomials
$U_{n,m}^{(k)}(z,w;a,b)$ which are the natural generalization of
the Umemura polynomials $U_n(z,w;a,b)$ related to the Painlev\'e
$VI$ equation. We show that if either $a=b$, or $a=0$, or $b=0$,
then polynomials $U_{n,m}^{(0)}(z,w;a,b)$ generate solutions to
the Painlev\'e $VI$ equation. We give new proof of
Noumi--Okada--Okamoto--Umemura conjecture, and describe
connections between polynomials $U_{n,m}^{(0)}(z,w;a,0)$ and
certain Umemura polynomials $U_k(z,w;\alpha,\beta)$. Finally we
show that after appropriate rescaling, Umemura's polynomials
$U_k(z,w;a,b)$ satisfy the Hirota--Miwa bilinear equations.
\end{abstract}
\section{Introduction}
There exists a vast body of literature about the Painlev\'e $VI$ equation
$P_{VI}:=P_{VI}(\alpha,\beta,\gamma,\delta):$
\begin{eqnarray}
\frac{d^2q}{dt^2}&=&\frac{1}{2}
\left(\frac{1}{q}+\frac{1}{q-1}+\frac{1}{q-t}\right)
\left(\frac{dq}{dt}\right)^2
-
\left(\frac{1}{t}+\frac{1}{t-1}+\frac{1}{q-t}\right)
\left(\frac{dq}{dt}\right) \nonumber \\
\smallskip \nonumber\\
&&+
\frac{q(q-1)(q-t)}{t^2(t-1)^2}
\left(
\alpha-\beta\frac{t}{q^2}
+\gamma\frac{(t-1)}{(q-1)^2}
+\delta\frac{t(t-1)}{(q-t)^2}
\right)
\end{eqnarray}
where $t \in {\bf C}$, $q:=q(t;\alpha,\beta,\gamma,\delta)$ is a
function of $t$, and $\alpha,\beta,\gamma,\delta$ are arbitrary
complex parameters, see e.g. \cite{NOOU,OI-OIV,P,U} and the
literature quoted therein. It is well--known and goes back to
Painlev\'e that any solution $q(t)$ of the equation $P_{VI}$
satisfies the so--called Painlev\'e property:
\begin{itemize}
\item
the critical points $0,1$ and $\infty$ of the equation (1.1) are
the only \emph{fixed singularities} of $q(t)$.
\item
any \emph{movable singularity} of $q(t)$
(the position of which depends on integration constants)
is a pole.
\end{itemize}
In this paper we introduce and initiate the study of certain
special polynomials related to the Painlev\'e $VI$ equation,
namely, the generalized Umemura polynomials
$U_{n,m}^{(k)}(z,w;a,b)$. These polynomials have many interesting
combinatorial and algebraic properties and in the particular case
$n=0=k$ coincide with Umemura's polynomials $U_m(z^2,w^2;a,b)$,
see e.g. \cite{U, NOOU}. The main goal of the present paper is to
study certain recurrence relations between polynomials
$U_{n,m}^{(k)}(z,w;a,b)$. Our main result is Theorem 1 which
gives a generalization of the recurrence relation between
Umemura's polynomials \cite{U}. In some particular cases the
recurrence relation obtained in Theorem 1 coincides with that for
Umemura's polynomials. As a corollary, we obtain a new proof of
the Noumi--Okada--Okamoto--Umemura conjecture \cite{NOOU}, and
show that polynomials $U_{n,m}^{(0)}(z,w;a,0)$ also generate
solutions to the equation Painlev\'e $VI$. The main mean in our
proofs is Lemma {\ref{b.lam.}} from Section 4. For example, using
this Lemma, we prove a new recurrence relation between Umemura's
polynomials (Theorem 2), describe explicitly connections between
polynomials $U_{n,m}(0,b)$ and Umemura's polynomials
$U_m(b_1,b_2)$, see Lemma {\ref{con.ume.}}, and prove that after
appropriate rescaling Umemura's polynomials $U_n(z,w;a,b)$
satisfy the Hirota--Miwa bilinear equations, see Proposition 5.
Finally, in Section 5, Proposition 6, we state and prove an
analog of the Pl\"ucker relations between certain Umemura's
polynomials.
\section{Painlev\'e $VI$}
\setcounter{equation}{0}
\renewcommand{\thesection}{\arabic{section}}
In this section we collect together some basic results about the
equation Painlev\'e $VI$. More detail and proofs may be found in
familiar series of papers by K.~Okamoto \cite{OI-OIV}. We refer
the reader to the Proceedings of Conference "The Painlev\'e
property. One century later" \cite{P}, where different aspects of
the theory of Painlev\'e equations may be found.
\subsection{Hamiltonian form}
It is well--known and goes back to a paper by Okamoto
\cite{OI-OIV} that the sixth Painlev\'e equation (1.1) is
equivalent to the following Hamiltonian system:
\begin{eqnarray}
\HH_{VI}( \bbb ;t,q,p):
\left\{
\begin{array}{lcr}
\frac{\displaystyle dq}{\displaystyle dt} &=&
\frac{\displaystyle \partial H}{\displaystyle \partial p},\\
\smallskip \\
\frac{\displaystyle dp}{\displaystyle dt} &=&
-\frac{\displaystyle \partial H}{\displaystyle \partial q},
\end{array}
\right.
\end{eqnarray}
with the Hamiltonian
\begin{eqnarray}
\lefteqn{H:=H_{VI}(\bbb;t,q,p)
=\frac{1}{t(t-1)}
\left[
q(q-1)(q-t)p^2
\right.
- \left\{
(b_1+b_2)(q-1)(q-t)
\right.
}
\nonumber\\
\smallskip \nonumber\\
&&
\left. \left.
+(b_1-b_2)q(q-t)
+(b_3+b_4)q(q-1)
\right\}p
+(b_1+b_3)(b_1+b_4)(q-t)
\right],\nonumber
\end{eqnarray}
where $\bbb=(b_1,b_2,b_3,b_4)$ belongs to the parameters space $\C^4$;
the parameters $(\alpha,\beta, \gamma, \delta)$ and $(b_1,b_2,b_3,b_4)$
are connected by the following relations
\begin{eqnarray}
\alpha=\frac{1}{2}(b_3-b_4)^2,
\beta=-\frac{1}{2}(b_1+b_2)^2,
\gamma=\frac{1}{2}(b_1-b_2)^2,
\delta=-\frac{1}{2}(b_3-b_4)(b_3+b_4-2).
\end{eqnarray}
\begin{Proposition}(K. Okamoto, \cite{OI-OIV}.)
If $(q(t),p(t))$ is a solution to the Hamiltonian system (2.1),
the function
\begin{eqnarray*}
h(\bbb,t)=t(t-1)H_{VI}(\bbb;t,q(t),p(t))
+e_2(b_1,b_3,b_4)t-\frac{1}{2}e_2(b_1,b_2,b_3,b_4)
\end{eqnarray*}
satisfies the equation $E_{VI}(\bbb):$
\begin{eqnarray}
\frac{dh}{dt} \left[t(t-1)\frac{d^2h}{dt^2}\right]^2
+\left[\frac{dh}{dt}\left\{2h-(2t-1)\frac{dh}{dt}
\right\}+b_1b_2b_3b_4\right]
=\prod_{k=1}^{4}\left(\frac{dh}{dt}+b_k^2\right),
\end{eqnarray}
where $e_2(x_1,\ldots ,x_n)=\displaystyle\sum_{1\le i<j\le
n}x_ix_j$ denotes the degree 2 elementary symmetric polynomial.
Conversely, for a solution $h:=h(\bbb,t)$ to the equation
$E_{VI}(\bbb)$ such that $\displaystyle\frac{d^2h}{dt^2} \neq 0$,
there exists a solution $(q(t),p(t))$ to the Hamiltonian system
(2.1). Furthermore, the function $q:=q(t)$ is a solution to the
Painlev\'e equation (1.1), where parameters $(\alpha, \beta,
\gamma, \delta)$ are determined by the relations (2.2).

\end{Proposition}
We will call the equation $E_{VI}(\bbb)$ by the
Painlev\'e--Okamoto equation.
\subsection{B\"acklund transformation}
Consider the following linear transformation of the parameters space $\C^4$:
\begin{eqnarray*}
&&s_1:=(b_1,b_2,b_3,b_4)\longmapsto (b_2,b_1,b_3,b_4),\\
\smallskip \\
&&s_2:=(b_1,b_2,b_3,b_4)\longmapsto (b_1,b_3,b_2,b_4),\\
\smallskip \\
&&s_3:=(b_1,b_2,b_3,b_4)\longmapsto (b_1,b_2,b_4,b_3),\\
\smallskip \\
&&s_0:=(b_1,b_2,b_3,b_4)\longmapsto (b_1,b_2,-b_3,-b_4),\\
\smallskip \\
&&l_3:=(b_1,b_2,b_3,b_4)\longmapsto (b_1,b_2,b_3+1,b_4).
\end{eqnarray*}
Denote by $W=<s_0,s_1,s_2,s_3,l_3>$ the subgroup of ${\rm Aut} \C^4$
generated by these transformation.
It is not difficult to see, that $W \cong W(D_4^{(1)})$, i.e.
$W$ is isomorphic to the affine Weyl group of type $D_4^{(1)}$.
\begin{Proposition}(K.~Okamoto, \cite{OI-OIV}.)
For each $w \in W$, there exists a birational transformation
\begin{eqnarray*}
L_w:\{\mbox{solutions to ${\cal H}_{VI}(\bbb)$ } \} \longmapsto
\{\mbox{solutions to ${\cal H}_{VI}(w(\bbb))$ } \}.
\end{eqnarray*}
\end{Proposition}
The birational transformations $L_w$, $w \in W(D_4^{(1)})$
are called by {\it B\"acklund transformations} associated to the equation
Painlev\'e $VI$.
\subsection{$\tau$--function}
Let $(q(t),p(t))$ be a solution to the Hamiltonian system (2.1),
the $\tau$--function $\tau(t)$ corresponding to the solution
$(q(t),p(t))$ is defined by the following equation
\begin{eqnarray*}
\frac{d}{dt} \log \tau(t) =H_{VI}(\bbb;t,q(t),p(t));
\end{eqnarray*}
in other words,
\begin{eqnarray*}
\tau(t)=({\rm constant}) \exp \left( \int H_{VI}(\bbb;t,q(t),p(t)) dt\right).
\end{eqnarray*}
\subsection{Umemura polynomials}
Suppose that $b_3=-\frac{1}{2}$, $b_4=0$, then it is well--known
and goes back to Umemura's paper \cite{U} , that the pair
\begin{eqnarray*}
(q_0,p_0)=\left(
\frac{(b_1+b_2)^2-(b_1^2-b_2^2)\sqrt{t(1-t)}}{(b_1-b_2)^2+4 b_1 b_2 t},
\frac{b_1q_0-\frac{1}{2}(b_1+b_2)}{q_0(q_0-1)}
\right)
\end{eqnarray*}
defines a solution to the Hamiltonian system (2.1) with parameters
$\bbb=(b_1,b_2,-\frac{1}{2},0)$. Note, see e.g. \cite{U}, that
\begin{eqnarray*}
\lefteqn{H_0(t)=H_{VI}\left((b_1,b_2,-\frac{1}{2},0);t,q_0(t),p_0(t)\right)}\\
\smallskip \\
&=&\frac{1}{t(t-1)}
\left\{
b_1(b_1-1)(1-2t)+2 b_1^2 \sqrt{t(t-1)}+2b_2(b_1-t) \right.\\
\smallskip \\
&&\left.
 +b_2(b_2-1)(1-2t)
-2b_2^2 \sqrt{t(t-1)}
\right\},
\end{eqnarray*}
and \begin{eqnarray*}
\tau_0(t)=\exp
\left\{\int H_0(t)dt
\right\}.
\end{eqnarray*}
To introduce Umemura's polynomials, let $(q_m,p_m)$ be a solution
to the Hamiltonian system ${\cal
H}_{VI}(b_1,b_2,-\frac{1}{2}+m,0)= {\cal
H}_{VI}(l_3^m(b_1,b_2,-\frac{1}{2},0)$ obtained from the solution
$(q_0,p_0)$ by applying $m$ times the the B\"acklund
transformation $l_3$. Consider the corresponding $\tau$--function
$\tau_m$:
\begin{eqnarray*}
\frac{d}{dt} \log \tau_m =H_{VI}((b_1,b_2,-\frac{1}{2}+m,0);t,q_m(t),p_m(t)).
\end{eqnarray*}
It follows from Proposition 1, see e.g. \cite{OI-OIV,U}, that
$\tau$--functions $\tau_n:=\tau_n(t)$ satisfy the Toda equation
\begin{eqnarray}
\frac{\tau_{n-1} \tau_{n+1}}{\tau_n^2}=
\frac{d}{dt}
\left( t(t-1)\frac{d}{dt}(\log \tau_n)
\right)+(b_1+b_2+n)(b_3+b_4+n).
\end{eqnarray}
Follow H. Umemura \cite{U}, define a family of functions $T_n(t)$,
$n=0,1,2,\ldots$, by
\begin{eqnarray*}
\tau_n(t)=T_n(t) \exp
\left( \int
\left( H_0(t)-\frac{n(b_1t-\frac{1}{2}(b_1+b_2))}{t(t-1)} \right)dt
\right).
\end{eqnarray*}
\begin{Proposition}(H. Umemura, \cite{U})
$T_n(t)$ is a polynomial with rational coefficients in the
variable
$v:=\displaystyle\sqrt{\frac{t}{t-1}}+\sqrt{\frac{t-1}{t}}$.
\end{Proposition}
For example, $T_0=1$, $T_1=1$,
$T_2=\frac{1}{2}\left( -4b_1^2+1)(2-v)/4+(-4b_2^2+1)(2+v)/4\right)$.\\
It follows from the Toda equation (2.4) that polynomials
$T_n:=T_n(v)$ satisfy the following recurrence relation \cite{U}:
\begin{eqnarray}
T_{n-1}T_{n+1} &=& \left \{
\frac{1}{4}(-2b_1^2-2b_2^2+(b_1^2-b_2^2)v)+(n-\frac{1}{2})^2
\right\}
T_n^2 \\
\smallskip \nonumber\\
&& + \frac{1}{4}(v^2-4)^2\left\{T_n \frac{d^2T_n}{dv^2}
-\left(\frac{d T_n}{dv}\right)^2\right\} \nonumber \\
\smallskip \nonumber\\
&& + \frac{1}{4}(v^2-4)v T_n \frac{d T_n}{dv}\nonumber
\end{eqnarray}
with initial conditions $T_0=T_1=1$.
\begin{Definition}
Polynomials $U_n:=U_n(z,w,b_1,b_2):=2^{n(n-1)}T_n(v)$, where
$z=\displaystyle\frac{2-v}{4}$,\break
$w=\displaystyle\frac{2+v}{4}$, are called by Umemura polynomials.
\end{Definition}
The formula (2.6) below was stated as a conjecture by M. Noumi,
S.~Okada,  K.~Okamoto and H.~Umemura \cite{NOOU} and has been
proved recently by M.~Taneda, and A.N.~Kirillov (independently):
\begin{eqnarray}
2^{n(n-1)}T_n(v) :=U_n(z,w,b_1,b_2)= \sum_{I \subset[n-1]}
\dd_n(I)
c_I d_{[n-1] \backslash I}
 \ z^{|I|} \  w^{|I^{c}|},
\end{eqnarray}
where
\begin{description}
\item[(i)] $[n-1]=\{1,2,\ldots,n-1\}$;
for any subset $I=\{i_1>i_2>\cdots>i_p\} \subset [n-1]$,\break
$\dd_n(I):=\dim_{\lambda(I)}^{GL(n)}$ stands for the dimension of
irreducible representation of the general linear group $GL(n)$
corresponding to the highest weight $\lambda(I)$ with the
Frobenius' symbol $\lambda(I)=(i_1,i_2, \ldots,i_p|i_1-1,i_2-1,
\ldots, i_p-1)$;
\item[(ii)]$c=-4b_1^2$, $d=-4b_2^2$, $z=\displaystyle\frac{2-v}{4}$,
$w=\displaystyle\frac{2+v}{4}$;
\item[(iii)]$\bar{c_k}=c+(2k-1)^2$, $\bar{d_k}=d+(2k-1)^2$,
$c_k=\bar{c}_1 \bar{c}_2 \cdots \bar{c}_k$,
$d_k=\bar{d}_1 \bar{d}_2 \cdots \bar{d}_k$;
\item[(iv)]
$|I|=i_1+i_2+\cdots+i_p$.
\end{description}
Recall that Frobenius' symbol $(a_1,a_2,\ldots,a_p|b_1,b_2,\ldots,b_p)$
denotes the partition which corresponds to the following diagram
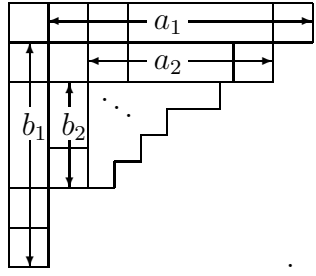
\begin{center}
\vskip2mm
\begin{picture}(120,100)
\put(10,100){\line(1,0){115}}
\put(10,85){\line(1,0){115}}
\put(10,70){\line(1,0){100}}
\put(25,45){\line(1,0){15}}
\put(10,30){\line(1,0){40}}
\put(10,15){\line(1,0){15}}
\put(10,0){\line(1,0){15}}

\put(50,30){\line(0,1){10}}
\put(50,40){\line(1,0){10}}
\put(60,40){\line(0,1){10}}
\put(60,50){\line(1,0){10}}
\put(70,50){\line(0,1){10}}
\put(70,60){\line(1,0){20}}
\put(90,60){\line(0,1){10}}

\put(10,100){\line(0,-1){100}}
\put(25,100){\line(0,-1){100}}
\put(40,100){\line(0,-1){70}}
\put(55,100){\line(0,-1){30}}
\put(110,100){\line(0,-1){15}}
\put(125,100){\line(0,-1){15}}
\put(95,85){\line(0,-1){15}}
\put(110,85){\line(0,-1){15}}

\put(62,93){\vector(-1,0){37}}
\put(78,93){\vector(1,0){47}}
\put(65,90){$a_1$}

\put(62,78){\vector(-1,0){22}}
\put(78,78){\vector(1,0){32}}
\put(65,75){$a_2$}

\put(18,61){\vector(0,1){24}}
\put(18,48){\vector(0,-1){48}}
\put(15,50){$b_1$}

\put(33,61){\vector(0,1){9}}
\put(33,48){\vector(0,-1){18}}
\put(30,50){$b_2$}

\put(44,56){$\ddots$}
\put(115,0){.}

\end{picture}
\end{center}
\renewcommand{\thesection}{\S \arabic{section}.}
\section{Generalized Umemura polynomials}
\setcounter{equation}{0} Let $n,m,k$ be fixed nonnegative
integers, $k \leq n$. Denote by $[n;m]$ the set of integers
$\{1,2,\ldots,n,n+2,n+4,\ldots,n+2m \}$. Let $I$ be a subset of
the set $[n;m]$. Follow \cite{DK}, define the numbers
\begin{eqnarray}
\dd_{n,m}(I)=\prod_{i \in I, j \in [n;m]\backslash I}
\left|\frac{i+j}{i-j}\right|,
\ \ \ \ \ \cc(I)=\sum_{i \in I, i>n} \frac{i-n}{2}
\end{eqnarray}
It has been shown in \cite{DK}, that in fact $\dd_{n,m}(I)$ are
integers for any subset $I \subset [n;m]$. Now we are going to
introduce the generalized Umemura polynomials
\begin{eqnarray*}
U_{n,m}^{(k)}:=U_{n,m}^{(k)}(z,w;a,b)=
\sum_{[k] \subset I \subset [n;m]}
\ \
\prod_{i \in I \backslash [k], j \in [k]}
\left(\frac{i+j}{i-j}\right)
\dd_{n,m}(I)
(-1)^{\cc(I)} e_{I}^{(n,m,k)}(z,w),
\end{eqnarray*}
where \\
\begin{description}
\item[(i)] $[k]$ stands for the set $\{1,2,\ldots,k\}$;
\item[(ii)]
$\bar{a}_k=a+(k-1)^2$, $\bar{b}_k=b+(k-1)^2$ and
$a_{2k}=\bar{a}_2 \bar{a}_4 \cdots \bar{a}_{2k}$,
$a_{2k+1}=\bar{a}_1 \bar{a}_3 \cdots \bar{a}_{2k+1}$;
$b_{2k}=\bar{b}_2 \bar{b}_4 \cdots \bar{b}_{2k}$,
$b_{2k+1}=\bar{b}_1 \bar{b}_3 \cdots \bar{b}_{2k+1}$;
\item[(iii)] for any subset $I \subset [n;m]$, we set
$a_I =\prod_{i \in I} a_i $, $b_I = \prod_{i \in I} b_i$;
\item[(iv)] $e_I^{(n,m,k)}(z,w)=a_{I \backslash [k]}
b_{[n;m]\backslash I}z^{|I\backslash [k]|}
w^{|[n;m]\backslash I|}$.
\end{description}
Note that the polynomial $U_{0,m}^{(0)}$ coincides with Umemura's
polynomial $T_m(z^2,w^2;a,b)$. The formula for generalized
Umemura polynomials stated below follows from the Cauchy
identity, and was used by
 J.F. van Diejen and A.N.~Kirillov \cite{DK}
in their study of $q$--spherical functions.
\begin{Lemma}{\label{det.rep.}}
The generalized Umemura polynomials $U_{n,m}^{(k)}(a,b;z,w)$
admit the following determinantal expression
$$U_{n,m}^{(k)}(a,b;z,w)= \det \left| a_i w^i \prod_{s \in
[k]}\left( \frac{i+s}{i-s} \right) \delta_{i,j}
+\frac{2i}{i+j}(-1)^{c({i})} \prod_{s \in [n;m], s \neq i} \left|
\frac{i+s}{i-s} \right| b_i z^i \right|_{i,j\in [n;m] \backslash
[k]},
$$
where $c(i)=i \ \ \mbox{if $i \leq n$, and }\ \  c(i)=(i-n)/2
\ \ \mbox{if $i >n$}$.
\end{Lemma}
In the particular case $k=0$, $n=0$ this formula gives a determinantal
representation for Umemura's polynomials and has many applications.
\section{Main result}
\setcounter{equation}{0}
Let us introduce notation $U_{n,m}:=U_{n,m}^{(0)}(z,w;a,b)$.
The main result of our paper describes a
recurrence relation between polynomials
$U_{n,m}$.
\begin{Theorem}
\begin{eqnarray}
U_{n,m-1}U_{n,m+1}&=&
\left( -\bar{a}_{n+2m+2}z^2+\bar{b}_{n+2m+2}w^2\right)U_{n,m}^2+
8z^2w^2D_x^2 U_{n,m} \circ U_{n,m} \nonumber \\
\smallskip \nonumber\\
&&
-\frac{4}{(n+2m+1)^2}ab(a-b)z^2w^2\left( U_{n,m}^{(1)}\right)^2,
\end{eqnarray}
where for any two functions $f=f(x)$ and $g=g(x)$
\begin{eqnarray*}
D_x^2f \circ g=f''g-2f'g'+fg''
\end{eqnarray*}
denotes the second Hirota derivative, and
$'=\displaystyle\frac{d}{dx}$; here variables $z$,$w$ and $x$ are
connected by the relations $z=\frac{1}{2}\left(
e^x+e^{-x}-2\right)^{1/2}$, $w=\frac{1}{2}\left(
e^x+e^{-x}+2\right)^{1/2}$.
\end{Theorem}
The main step of our proof is to establish the following algebraic
identity which appears to have an independent interest.
\begin{Lemma}{\label{b.lam.}}
For any two subsets $I$, $J$ of the set $[n;m]$, we have
\begin{eqnarray*}
&&\prod_{ \lambda \in I}
\left(
\frac{ x + 2+\lambda}{ x +2- \lambda}
\right)
\prod_{ \lambda \in J}
\left(
\frac{ x - \lambda}{x  + \lambda}
\right)+
\prod_{ \lambda \in I}
\left(
\frac{ x - \lambda }{ x + \lambda }
\right)
\prod_{ \lambda \in J}
\left(
\frac{ x +2+ \lambda }{ x +2- \lambda }
\right)
\\
\smallskip \nonumber\\
&&\ \ \ = 2 + \sum_{ \lambda \in I \cup J}
\frac{b_{ \lambda }^{ I , J }}{( x +2- \lambda )( x + \lambda )},\nonumber
\end{eqnarray*}
where the coefficients $b_{\lambda}^{I,J}$ have the following expressions:
\begin{description}
\item[(i)] if $\lambda \neq 1$ and $\lambda \in I \cap J$, then
$b_\lambda^{I,J}=$
$$\!\!\!\!\!\!\!\!\!\!\! 4 \lambda (\lambda-1) \left\{ \prod_{
\lambda' \in I \backslash \{ \lambda \}} \left(\frac{ \lambda
+\lambda'}{ \lambda - \lambda'} \right) \prod_{ \lambda' \in J}
\left( \frac{ \lambda - 2 - \lambda'}{\lambda - 2  + \lambda'}
\right)+ \prod_{ \lambda' \in I} \left( \frac{ \lambda - 2 -
\lambda' }{ \lambda - 2 + \lambda' } \right) \prod_{ \lambda' \in
J \backslash \{\lambda \} } \left( \frac{ \lambda + \lambda' }{
\lambda - \lambda' } \right) \right\};
$$
\item[(ii)] if $\lambda \neq 1$, $\lambda \in I$ and $\lambda \not\in J$,
then
\begin{eqnarray*}
b_\lambda^{I,J}=
4 \lambda (\lambda-1) \left\{
\prod_{ \lambda' \in I \backslash \{ \lambda \}}
\left(
\frac{ \lambda +\lambda'}{ \lambda - \lambda'}
\right)
\prod_{ \lambda' \in J}
\left(
\frac{ \lambda - 2 - \lambda'}{\lambda - 2  + \lambda'}
\right)
\right\};
\end{eqnarray*}
\item[(iii)] if $1 \in I \cap J$, then
\begin{eqnarray*}
b_1^{I,J}=
-8 \prod_{ \lambda \in I \backslash \{ 1 \}}
\left(
\frac{ 1 +\lambda}{ 1 - \lambda}
\right)
\prod_{ \lambda \in J \backslash \{ 1 \} }
\left(
\frac{ 1+ \lambda}{1  - \lambda}
\right).
\end{eqnarray*}
\end{description}
\end{Lemma}
{\bf Proof.}  Using the partial fraction expansion, we have
\begin{eqnarray*}
&&\prod_{ \lambda \in I}
\left(
\frac{ x + 2+\lambda}{ x +2- \lambda}
\right)
\prod_{ \lambda \in J}
\left(
\frac{ x - \lambda}{x  + \lambda}
\right)+
\prod_{ \lambda \in I}
\left(
\frac{ x - \lambda }{ x + \lambda }
\right)
\prod_{ \lambda \in J}
\left(
\frac{ x +2+ \lambda }{ x +2- \lambda }
\right)
\\
\smallskip \\
&=&
2+\frac{A}{x+1}+\frac{B}{(x+1)^2}+
\sum_{\lambda \in I,\lambda \neq 1} \frac{C_\lambda}{x+2+\lambda}\\
\smallskip \\
&&+\sum_{\lambda \in I,\lambda \neq 1} \frac{D_\lambda}{x-\lambda}
+\sum_{\lambda \in J,\lambda \neq 1} \frac{E_\lambda}{x+2+\lambda}
+\sum_{\lambda \in J,\lambda \neq 1} \frac{F_\lambda}{x-\lambda},
\end{eqnarray*}
where $A$, $B$, $C_\lambda$, $D_\lambda$, $E_\lambda$ and $F_\lambda$
are some constants. It follows from the residue theorem that we have
\begin{eqnarray*}
A=0 \mbox{ and } C_\lambda=-D_\lambda=
2\lambda
\left\{
\prod_{ \lambda' \in I \backslash \{ \lambda \}}
\left(
\frac{ \lambda +\lambda'}{ \lambda - \lambda'}
\right)
\prod_{ \lambda' \in J}
\left(
\frac{ \lambda - 2 - \lambda'}{\lambda - 2  + \lambda'}
\right)
\right\}.
\end{eqnarray*}
Similarly, we get expressions for $E_\lambda$ and $F_\lambda$.
Moreover, If $1 \in I \cap J$, then
\begin{eqnarray*}
B=-8 \prod_{ \lambda \in I \backslash \{ 1 \}}
\left(
\frac{ 1 +\lambda}{ 1 - \lambda}
\right)
\prod_{ \lambda \in J \backslash \{ 1 \} }
\left(
\frac{ 1+ \lambda}{1  - \lambda}
\right).
\end{eqnarray*}
All the statements of Lemma {\ref{b.lam.}} follow from the above expressions
for coefficients $C_\lambda$, $D_\lambda$ and $B$ by direct calculations.
\begin{Lemma}{\label{sum.b.lam.}}
For any two subsets $I$, $J$ of the set $[n;m]$, we have
\begin{eqnarray}{\label{sum.b.lam.equ.}}
\sum_{ \lambda \in I \cup J }
b_{ \lambda }^{ I , J } = 4( |I| - |J| )^2 - 4 (|I| + |J|).
\end{eqnarray}
\end{Lemma}
{\bf Proof.} By Lemma {\ref{b.lam.}}, we have
\begin{eqnarray*}
\sum_{ \lambda \in I \cup J }
b_{ \lambda }^{ I , J } &=&\lim_{x \rightarrow \infty}
\left\{ (\mbox{L.H.S of Lemma {\ref{sum.b.lam.}}
{(\ref{sum.b.lam.equ.})} })-2 \right\}x^2\\
 \smallskip \\
&=& \lim_{x \rightarrow \infty}
x^2\left\{
\prod_{ \lambda \in I}
\left(
\frac{ x + 2+\lambda}{ x +2- \lambda}
\right)
\prod_{ \lambda \in J}
\left(
\frac{ x - \lambda}{x  + \lambda}
\right)+
\prod_{ \lambda \in I}
\left(
\frac{ x - \lambda }{ x + \lambda }
\right)
\prod_{ \lambda \in J}
\left(
\frac{ x +2+ \lambda }{ x +2- \lambda }
\right)
-2\right\}\\
\smallskip \\
&=&\lim_{x \rightarrow \infty}
\frac{\displaystyle A(x)}
{\displaystyle
\prod_{\lambda \in I} (x+2-\lambda)
\prod_{\lambda \in I} (x+\lambda)
\prod_{\lambda \in J} (x+2-\lambda)
\prod_{\lambda \in J} (x+\lambda)
},
\end{eqnarray*}
where
\begin{eqnarray*}
A(x)&=&x^2 \left\{\prod_{\lambda \in I} (x+2+\lambda)
\prod_{\lambda \in I} (x+\lambda)
\prod_{\lambda \in J} (x+2-\lambda)
\prod_{\lambda \in J} (x-\lambda)\right.\\
\smallskip \\
&&+\prod_{\lambda \in I} (x+2-\lambda)
\prod_{\lambda \in I} (x-\lambda)
\prod_{\lambda \in J} (x+2+\lambda)
\prod_{\lambda \in J} (x+\lambda)\\
\smallskip \\
&&\left.-2 \prod_{\lambda \in I} (x+2-\lambda)
\prod_{\lambda \in I} (x+\lambda)
\prod_{\lambda \in J} (x+2-\lambda)
\prod_{\lambda \in J} (x+\lambda)\right\}.
\end{eqnarray*}
It is easy to see that
coefficients of $x^{|I|+|J|+2}$ and $x^{|I|+|J|+1}$ in $A(x)$ are disappear.
Now let us compute the coefficient of $x^{|I|+|J|}$ in $A(x)$:
\begin{eqnarray*}
&&\sum_{\lambda \in I}
\left\{(2+\lambda)(\lambda)+(2-\lambda)(-\lambda)-2(2-\lambda)(\lambda)
\right\}\\
\smallskip \nonumber\\
&&+\sum_{\lambda_1, \lambda_2 \in I, \lambda_1 <\lambda_2}
\left\{
(2+\lambda_1)(2+\lambda_2)
+(2+\lambda_1)(\lambda_2)
+(2+\lambda_2)(\lambda_1)
+(\lambda_1)(\lambda_2)\right.\\
\smallskip \nonumber
&&+(2-\lambda_1)(2-\lambda_2)
+(2-\lambda_1)(-\lambda_2)
+(2-\lambda_2)(-\lambda_1)
+(-\lambda_1)(-\lambda_2)\\
\smallskip \\
&&
\left.-2(2-\lambda_1)(2-\lambda_2)
-2(2-\lambda_1)(\lambda_2)
-2(2-\lambda_2)(\lambda_1)
-2(\lambda_1)(\lambda_2)
\right\}\\
\smallskip \\
&&+\sum_{\lambda_1 \in I, \lambda_2 \in J}
\{(2+2\lambda_1)(2-2\lambda_2)+(2-2\lambda_1)(2+2\lambda_2)-8
\}\\
\smallskip \\
&&+\sum_{\lambda \in J}
\left\{(2+\lambda)(\lambda)+(2-\lambda)(-\lambda)-2(2-\lambda)(\lambda)
\right\}\\
\smallskip \nonumber\\
&&+\sum_{\lambda_1, \lambda_2 \in J, \lambda_1 <\lambda_2}
\left\{
(2+\lambda_1)(2+\lambda_2)
+(2+\lambda_1)(\lambda_2)
+(2+\lambda_2)(\lambda_1)
+(\lambda_1)(\lambda_2)\right.\\
\smallskip \nonumber
&&+(2-\lambda_1)(2-\lambda_2)
+(2-\lambda_1)(-\lambda_2)
+(2-\lambda_2)(-\lambda_1)
+(-\lambda_1)(-\lambda_2)\\
\smallskip \\
&&
\left.-2(2-\lambda_1)(2-\lambda_2)
-2(2-\lambda_1)(\lambda_2)
-2(2-\lambda_2)(\lambda_1)
-2(\lambda_1)(\lambda_2)
\right\}\\
\smallskip \\
&=&\sum_{\lambda \in I} (4 \lambda^2 -4 \lambda)
+\sum_{\lambda_1,\lambda_2 \in I, \lambda_1< \lambda_2} 8 \lambda_1 \lambda_2
-\sum_{\lambda_1 \in I, \lambda_2 \in J}8 \lambda_1 \lambda_2\\
\smallskip \\
&&+\sum_{\lambda \in J} (4 \lambda^2 -4 \lambda)
+\sum_{\lambda_1,\lambda_2 \in J, \lambda_1< \lambda_2} 8 \lambda_1 \lambda_2
\\
\smallskip \\
&=&\left\{ 4(|I|-|J|)^2-4(|I|+|J|)\right\}.\\
\end{eqnarray*}
The latter expression coincide with the RHS (4.2), and therefore
the proof of Lemma 3 is finished.

\begin{Lemma}{\label{b.lam.0.}}
For an element $ \lambda \in I \cap J $, we have
$ b_{ \lambda } ^{ I , J }=0 $ if and only if $ \lambda - 2 \in I \cap J $.
For an element  $ \lambda \in I \backslash  (I \cap J) $, we have
$b_{ \lambda}^{I,J}=0 $ if and only if $ \lambda-2 \in J $.
\end{Lemma}
This lemma follows from Lemma {\ref{b.lam.}} by direct calculation.

\begin{Lemma}{\label{hir.dif.}}
If $ n_1+m_1=n_2+m_2 $, then
\begin{eqnarray}
\lefteqn{4 z^2 w^2 D_x^2 z^{n_1}w^{m_1} \circ z^{n_2} w^{m_2}
=
\left[ - \left\{ (n_1+n_2)-(n_1-n_2)^2 \right\} \right. w^2}
\\
\smallskip \nonumber\\
&&
+ \left. \left\{ (m_1+m_2)-(m_1-m_2)^2 \right\} z^2 \right]
z^{n_1+n_2} w^{m_1+m_2}. \nonumber
\end{eqnarray}
\end{Lemma}
Easy proof by direct computation.

Now we are ready to prove our main theorem.

Since $D_x^2$ is a bilinear operator, using the above identity
(4.3), we have
\begin{eqnarray*}
\lefteqn{\left( -\bar{a}_{n+2m+2}z^2+\bar{b}_{n+2m+2}w^2\right)U_{n,m}^2+
8z^2w^2D_x^2 U_{n,m} \circ U_{n,m} }
\\
\smallskip \nonumber\\
&=&\sum_{I,J \subset [n;m]}
P_b(I,J)w^2 \dd_{n,m}(I) \dd_{n,m}(J)(-1)^{\cc(I)+\cc(J)}
e_I^{(n,m,0)}(z,w) e_J^{(n,m,0)}(z,w)\\
\smallskip\\
&&-
\sum_{I,J \subset [n;m]}
P_a(I,J)z^2 \dd_{n,m}(I) \dd_{n,m}(J)(-1)^{\cc(I)+\cc(J)}
e_{[n;m]\backslash I}^{(n,m,0)}(z,w) e_{[n;m]\backslash J}^{(n,m,0)}(z,w),\\
\end{eqnarray*}
where we write
\begin{eqnarray*}
P_a(I,J)=\bar{a}_{n+2m+2}-2\left\{ \left(|I|+|J|\right)
-\left(|I|-|J| \right)^2\right\},\\
\smallskip \\
P_b(I,J)=\bar{b}_{n+2m+2}-2\left\{ \left(|I|+|J|\right)
-\left(|I|-|J| \right)^2\right\}.
\end{eqnarray*}
By Lemma {\ref{sum.b.lam.}} we have
\begin{eqnarray*}
P_a(I,J)&=&
\bar{a}_{n+2m+2}+\frac{1}{2}\sum_{\lambda \in {I \cup J}} b_{\lambda}^{I,J}\\
\smallskip \\
&=&\bar{a}_{n+2m+2}\left(1+\sum_{\lambda \in I \cup J }
\frac{b_\lambda^{I,J}}{2(n+2m+2-\lambda)(n+2m+\lambda)}\right)\\
\smallskip\\
&&-\sum_{\lambda \in {I \cup J}}
\frac{\bar{a}_{n+2m+2}-(n+2m+2-\lambda)(n+2m+\lambda)}
{2(n+2m+2-\lambda)(n+2m+\lambda)} b_\lambda^{I,J}.\\ && ~
\end{eqnarray*}
Using Lemma {\ref{b.lam.}}, the latter expression can be
transformed to the following form:
\begin{eqnarray*}
P_a(I,J)&=&
\frac{\bar{a}_{n+2m+2}}{2}
\left(
\prod_{ \lambda \in I}
\left(
\frac{ n+2m + 2+\lambda}{ n+2m +2- \lambda}
\right)
\prod_{ \lambda \in J}
\left(
\frac{ n+2m - \lambda}{n+2m  + \lambda}
\right)\right.\\
\smallskip\\
&&+\left.
\prod_{ \lambda \in I}
\left(
\frac{ n+2m - \lambda }{ n+2m + \lambda }
\right)
\prod_{ \lambda \in J}
\left(
\frac{ n+2m +2+ \lambda }{ n+2m +2- \lambda }
\right)
\right)\\
\smallskip\\
&&+\sum_{\lambda \in {I \cup J}}
\frac{\bar{a}_{\lambda}}
{2(n+2m+2-\lambda)(n+2m+\lambda)} b_\lambda^{I,J}.
\end{eqnarray*}
Now if $n+2m \in J$, then it is not difficult to check that
\begin{eqnarray*}
\lefteqn{\frac{\bar{b}_{n+2m+2}}{2}
\prod_{ \lambda \in I}
\left(
\frac{ n+2m + 2+\lambda}{ n+2m +2- \lambda}
\right)
\prod_{ \lambda \in J}
\left(
\frac{ n+2m - \lambda}{n+2m  + \lambda}\right)}\\
\smallskip \\
&&\times
w^2 \dd_{n,m}(I) \dd_{n,m}(J)(-1)^{\cc(I)+\cc(J)}
e_I^{(n,m,0)}(z,w) e_J^{(n,m,0)}(z,w)\\
\smallskip \\
&=&\frac{1}{2}\dd_{n,m+1}(I) \dd_{n,m-1}(J)(-1)^{\cc(I)+\cc(J)}
e_I^{(n,m+1,0)}(z,w) e_J^{(n,m-1,0)}(z,w).
\end{eqnarray*}
Hence we have
\begin{eqnarray*}
\lefteqn{\left( -\bar{a}_{n+2m+2}z^2+\bar{b}_{n+2m+2}w^2\right)U_{n,m}^2+
8z^2w^2D_x^2 U_{n,m} \circ U_{n,m} }
\\
\smallskip \nonumber\\
&=&\sum_{I \subset [n,m+1], \ J \subset [n,m-1]}
\dd_{n,m+1}(I) \dd_{n,m-1}(J)(-1)^{\cc(I)+\cc(J)}
e_I^{(n,m+1,0)}(z,w) e_J^{(n,m-1,0)}(z,w)\\
\smallskip\\
&&
+\sum_{I,J \in [n;m]} \sum_{\lambda \in {I \cup J}}
\frac{ \ b_\lambda^{I,J}
\dd_{n,m}(I) \dd_{n,m}(J)(-1)^{\cc(I)+\cc(J)}}
{2(n+2m+2-\lambda)(n+2m+\lambda)}
\bar{b}_{\lambda} w^2 e_I^{(n,m,0)}(z,w) e_J^{(n,m,0)}(z,w)
\\
\smallskip\\
&&-
\sum_{I \subset [n,m+1], \ J \subset [n,m-1]}
\dd_{n,m+1}(I) \dd_{n,m-1}(J)(-1)^{\cc(I)+\cc(J)}
e_{[n;m]\backslash I}^{(n,m+1,0)}(z,w) e_{[n;m] \backslash J}^{(n,m-1,0)}(z,w)
\\
\smallskip \\
&&-\sum_{I,J \in [n;m]} \sum_{\lambda \in {I \cup J}}
\frac{ \ b_\lambda^{I,J}
\dd_{n,m}(I) \dd_{n,m}(J)(-1)^{\cc(I)+\cc(J)}}
{2(n+2m+2-\lambda)(n+2m+\lambda)}
\bar{a}_{\lambda} z^2 e_{[n;m] \backslash I}^{(n,m,0)}(z,w)
e_{[n;m] \backslash J}^{(n,m,0)}(z,w)\\
\smallskip \\
&=&U_{n,m+1} U_{n,m-1}\\
\smallskip \\
&&+\sum_{I,J \in [n;m]} \sum_{\lambda \in {I \cup J}}
\frac{ \ b_\lambda^{I,J}
\dd_{n,m}(I) \dd_{n,m}(J)(-1)^{\cc(I)+\cc(J)}}
{2(n+2m+2-\lambda)(n+2m+\lambda)}
\bar{b}_{\lambda} w^2 e_I^{(n,m,0)}(z,w) e_J^{(n,m,0)}(z,w)
\\
\smallskip \\
&&-\sum_{I,J \in [n;m]} \sum_{\lambda \in {I \cup J}}
\frac{ \ b_\lambda^{I,J}
\dd_{n,m}(I) \dd_{n,m}(J)(-1)^{\cc(I)+\cc(J)}}
{2(n+2m+2-\lambda)(n+2m+\lambda)}
\bar{a}_{\lambda} z^2 e_{[n;m] \backslash I}^{(n,m,0)}(z,w)
e_{[n;m] \backslash J}^{(n,m,0)}(z,w).
\end{eqnarray*}
To continue our proof, for $\lambda \neq 1$, $\lambda \in I$, we
define the function $\Split(b_\lambda)(I,J)$ by
\begin{eqnarray*}
\Split(b_\lambda)(I,J)=4\lambda (\lambda-1)
\prod_{ \lambda' \in I \backslash \{ \lambda \}}
\left(
\frac{ \lambda +\lambda'}{ \lambda - \lambda'}
\right)
\prod_{ \lambda' \in J}
\left(
\frac{ \lambda - 2 - \lambda'}{\lambda - 2  + \lambda'}
\right).
\end{eqnarray*}
Let us observe that if $\lambda \in I \cap J$ and $\lambda \neq 1$, then
\begin{eqnarray*}
b_{\lambda}^{I,J}=\Split(b_{\lambda})(I,J)+\Split(b_{\lambda})(J,I).
\end{eqnarray*}
Similarly, if $\lambda \in I$, $\lambda \neq 1$ and $\lambda \not\in J$, then
$b_{\lambda}^{I,J}=\Split(b_{\lambda})(I,J)$.

Finally, if $\lambda \in I $, $\lambda \neq 1$ and $\lambda-2 \not\in J$, then
\begin{eqnarray*}
\prod_{\lambda' \in I \backslash \{ \lambda\}}
\left(
\frac{ \lambda +\lambda'}{ \lambda - \lambda'}
\right)
\dd_{n,m}(I)
=(-1)^{\#[n;m]_{> \lambda}}
\prod_{\lambda \in [n;m] \backslash \{I \backslash \{ \lambda\}\}}
\left(
\frac{\lambda+\lambda'}{\lambda-\lambda'}
\right)
\dd_{n,m}([n;m] \backslash \{I \backslash \{\lambda \}\})
\end{eqnarray*}
and
\begin{eqnarray*}
\lefteqn{
\prod_{\lambda' \in I \backslash \{ \lambda\}}
\left(
\frac{ \lambda -2-\lambda'}{ \lambda -2+ \lambda'}
\right)
\dd_{n,m}(I)
=(-1)^{\#[n;m]_{> \lambda-2}}}\\
\smallskip \\
&&\times
\prod_{\lambda \in [n;m] \backslash \{I \cup \{ \lambda-2\}\}}
\left(
\frac{\lambda-2-\lambda'}{\lambda-2+\lambda'}
\right)
\dd_{n,m}([n;m] \backslash \{I \cup \{\lambda-2 \}\}),
\end{eqnarray*}
where $\#[n;m]_{>\lambda}=\#\{i \in [n;m]|i > \lambda\}$.

Let us summarize the results of above calculations as an auxiliary lemma.
\begin{Lemma}{\label{dis.}}
If $\lambda \in I $, $\lambda \neq 1$ and $\lambda-2 \not\in J$, then
\begin{eqnarray*}
\Split(b_{\lambda})(I,J)\dd_{n,m}(I)\dd_{n,m}(J)
=(-1)^{A}
\Split(b_{\lambda})(I',J')
\dd_{n,m}(I')\dd_{n,m}(J'),
\end{eqnarray*}
where $I'=[n;m] \backslash \{I\backslash \{\lambda\}\} $,
$J'=[n;m] \backslash \{J \cup \{\lambda-2 \}\} $ and
$A=1$ if $\lambda \leq n$, and $A=-1$ if $\lambda >n$.
\end{Lemma}
Note that under the assumption of Lemma {\ref{dis.}}, we have
$\bar{b}_\lambda w^2 e_{I}^{(n,m,0)} e_{J}^{(n,m,0)}
=\bar{a}_\lambda z^2 e_{I'}^{(n,m,0)} e_{J'}^{(n,m,0)}$
and
$\bar{a}_\lambda w^2 e_{I}^{(n,m,0)} e_{J}^{(n,m,0)}
=\bar{b}_\lambda z^2 e_{I'}^{(n,m,0)} e_{J'}^{(n,m,0)}$.
Hence, using Lemma {\ref{dis.}}, we obtain the following equality
\begin{eqnarray*}
\lefteqn{\left( -\bar{a}_{n+2m+2}z^2+\bar{b}_{n+2m+2}w^2\right)U_{n,m}^2+
8z^2w^2D_x^2 U_{n,m} \circ U_{n,m} }
\\
\smallskip \nonumber\\
&=&U_{n,m+1} U_{n,m-1}\\
\smallskip \\
&&+\sum_{I,J \in [n;m],\ 1 \in I \cap J} \frac{ \ b_1^{I,J}
\dd_{n,m}(I) \dd_{n,m}(J)(-1)^{\cc(I)+\cc(J)}}
{2(n+2m+1)^2}
\bar{b}_{1} w^2 e_I^{(n,m,0)}(z,w) e_J^{(n,m,0)}(z,w)
\\
\smallskip \\
&&-\sum_{I,J \in [n;m], \ 1 \in I \cap J} \frac{ \ b_1^{I,J}
\dd_{n,m}(I) \dd_{n,m}(J)(-1)^{\cc(I)+\cc(J)}}
{2(n+2m+1)^2}
\bar{a}_{1} z^2 e_{[n;m] \backslash I}^{(n,m,0)}(z,w)
e_{[n;m] \backslash J}^{(n,m,0)}(z,w).
\end{eqnarray*}
To simplify the right hand side of the latter equality we are going
to use the following expression for the coefficients $b_1^{I,J}$
which is an easy consequence of Lemma 2:
\begin{eqnarray*}
b_1^{I,J}=
-8 \prod_{ \lambda' \in I \backslash \{ 1 \}}
\left(
\frac{ 1 +\lambda'}{ 1 - \lambda'}
\right)
\prod_{ \lambda' \in J \backslash \{ 1 \} }
\left(
\frac{ 1+ \lambda'}{1  - \lambda'}
\right).
\end{eqnarray*}
After substituting the above expressions for the coefficients
$b_1^{I,J}$ to the both sums which appear in the right hand side
of the equality under consideration, it remains to observe that if
$1 \in I \subset [n;m]$, then
\begin{eqnarray*}
\prod_{\lambda ' \in I \backslash \{1\}}
\left(\frac{1+\lambda'}{1-\lambda'}\right)
\dd_{n,m}(I)
=
\prod_{\lambda ' \in [n;m] \backslash I}
\left(\frac{1+\lambda'}{1-\lambda'}\right)
\dd_{n,m}(I \backslash \{1\}),
\end{eqnarray*}
\begin{eqnarray*}
a^2 b w^2 z^2 e_I^{(n,m,1)}(z,w) e_J^{(n,m,1)}(z,w)
=
\bar{b}_{1} w^2 e_I^{(n,m,0)}(z,w) e_J^{(n,m,0)}(z,w)
\end{eqnarray*} and
\begin{eqnarray*}
a b^2 w^2 z^2 e_I^{(n,m,1)}(z,w) e_J^{(n,m,1)}(z,w)
=
\bar{a}_{1} z^2 e_{[n;m] \backslash I}^{(n,m,0)}(z,w)
e_{[n;m] \backslash J}^{(n,m,0)}(z,w).
\end{eqnarray*}
The proof of Theorem 1 is finished.

\hfill$\square$

\begin{description}
\item[Remarks 1.]If $n=0$, then $U_{0,m}=U_{m+1}(z^2,w^2;a,b)$
coincides with the Umemura polynomial, and $U_{0,m}^{(1)}=0$. In
this case the recurrence relation (4.1) has been used by M.~Taneda
in his proof of Noumi--Okada--Okamoto--Umemura's Conjecture (2.6).
\item[2.]
Note that $U_{0,m}=U_{2,m-1}^{(1)}/(2m+1)$, and more generally
$$U_{k,m}^{(k)}=U_{k+2,m-1}^{(k+1)}(2k+1)!!(2m-1)!!/(2k+2m+1)!!,
$$
where $(2n+1)!!=1 \cdot 3 \cdot 5 \cdots (2n+1)$.
\item[3.]
"Unwanted term" in (4.1) which contains $\left( U_{n,m}^{(1)} \right)^2$
vanishes if either $a=0$, or $b=0$, or $a=b$.
\end{description}
In the case $a=b$ and $k=0$ the expression $e_I^{(n,m,k)}(z,w)$
doesn't depend on a subset $I \subset [n;m]$ and is equal to
$a_{[n;m]} z^{|I|}w^{|[n;m] \backslash I|}$. Hence,
in this case
\begin{eqnarray}
U_{n;m}^{(0)}(z,w;a,a)&=&
a_{[n;m]} \sum_{I \subset [n;m]} \dd_{n,m}(I) (-1)^{\cc(I)}
z^{|I|}w^{|[n;m] \backslash I|} \nonumber \\
\smallskip \nonumber \\
&=&a_{[n;m]}(z+w)^{\tiny  \left(
\begin{array}{c} n+m+1 \\2
\end{array}\right)}
(z-w)^{\tiny \left( \begin{array}{c} m+1 \\2
\end{array}\right) }.
\end{eqnarray}
Recall that $\bar{a}_i=a+(i-1)^2$, $a_{2i}=\bar{a}_2 \bar{a}_4
\cdots \bar{a}_{2i}$, $a_{2i+1}=\bar{a}_1 \bar{a}_3 \cdots
\bar{a}_{2i+1}$ and $a_{[n;m]}=\prod_{ i \in [n;m] } a_i$. The
last equality in (4.4) has been proved for the first time by
J.F.~van Diejen and A.N.~Ki\-rillov \cite{DK}. On the other hand,
we can show that polynomials
\begin{eqnarray*}
X_{n,m}(z,w;a)=a_{[n;m]}(z+w)^{\tiny  \left(\!\!
\begin{array}{c} n+m+1  \\2
\end{array}\!\!\right)}
(z-w)^{\tiny \left( \!\!\begin{array}{c} m+1 \\2
\end{array}\!\!\right) }
\end{eqnarray*}
also satisfy the recurrence relation (4.1) and coincide with
polynomials $U_{n,m}^{(0)}(z,w;a,a)$ if $m=0$. From this
observation we can deduce the equality
$X_{n,m}(z,w;a)=U_{n,m}^{(0)}(z,w;a,a)$, which is equivalent to
the main identity from \cite{DK}. Another case when "unwanted
term" in (4.1) vanishes is the case when either $a=0$, or $b=0$.
In this case we have
\begin{Corollary}
Assume that $b=0$, then polynomial $U_{n,m}(z,w;a,0)$
defines a solution to the equation Painlev\'e $VI$.
\end{Corollary}
Finally, we are going to compare polynomials $U_{n,m}(z,w;a,0)$ and
$U_m(z,w;\alpha,\beta)$.
For this goal, let us consider functions
\begin{eqnarray*}
h_0:=h_0(t)=\left\{ b_1^2\left( \sqrt{t}-\sqrt{t-1} \right)^2
+b_2^2\left( \sqrt{t}+\sqrt{t+1}\right)^2
\right\}/4,
\end{eqnarray*}
and
\begin{eqnarray*}
h_{n,m}:=h_{n,m}(b_1,b_2)=t(t-1) \log \left(U_{n,m} \right)'-h_0.
\end{eqnarray*}
\begin{Proposition}{\label{aux.ham.}}
\begin{description}
\item[(i)]  $h_{0,m}$ satisfies the Painlev\'e--Okamoto
equation $E_{VI}(b_1,b_2,m+1/2,0)$;
\item[(ii)] $h_{1,m}=-(2t-1)(m+1)^2/2$ satisfies
the equation $E_{VI}(0,m+1,b_3,b_4)$;
\item[(iii)] $h_{n,m}(0,b_2)$ satisfies the equation
$E_{VI}(0,b_2,\displaystyle\frac{n}{2},\frac{n+2m+1}{2})$.
\end{description}
\end{Proposition}
Proposition 4 follows from
Lemma {\ref{con.ume.}} and Lemma \ref{con.aux.ham.} below.

Let us define $U_{n,m}(b_1,b_2):=U_{n,m}(z,w;-4b_1^2,-4b_2^2)$, then
\begin{Lemma}{\label{con.ume.}}
\begin{eqnarray*}
U_{n,m}(0,b_2)
=\left\{
\begin{array}{lcr}
\displaystyle
b_{[n;m]_\odd} w^{\left(\frac{n}{2}\right)^2}
U_{0,m+\frac{n}{2}}\left(\frac{n}{2},b_2\right), & \mbox{if $n$ is even},
&(\theequation)\\
\smallskip \\ \addtocounter{equation}{1}
\displaystyle
b_{[n;m]_\odd} w^{\left(\frac{n+2m+1}{2}\right)^2}
U_{0,\frac{n-1}{2}}\left(m+\frac{n+1}{2},b_2\right), & \mbox{if $n$ is odd},
&(\theequation) \addtocounter{equation}{1}
\end{array}
\right.
\end{eqnarray*}
where $[n;m]_\odd=\left\{ i \in [n;m] \ \ |\ \  \mbox{$i$ is odd }\right\}$.
\end{Lemma}

{\bf Proof}. By the definition of generalized Umemura polynomials,
one can see
\begin{eqnarray*}
U_{n,m}(0,b_2)=
\sum_{I \subset [n;m]_\even}\dd_{n,m}{I} (-1)^{\cc(I)}e_I^{(n,m,0)}(z,w).
\end{eqnarray*}
Assume first that both $n$ and $i$ are even, then we have
\begin{eqnarray*}
\left.
a_i\right|_{b_1=0}=\{(i-1)!!\}^2,
\end{eqnarray*}
\begin{eqnarray*}
\left. a_i \right|_{b_1=n/2}=\prod_{j=1}^{i/2}\{-n^2+(2j-1)^2\}
\left. a_i\right|_{b_1=0}.
\end{eqnarray*}
Now if $i \leq n$, then
\begin{eqnarray*}
\prod_{j \in [n;m]_\odd}\left|\frac{i+j}{i-j} \right|
\left. a_i\right|_{b_1=0}
&=&
\frac{(i+n-1)!!}{(n-1-i)!!}
=(n-i+1)(n-i+3) \cdots (n-i-1)\\
\smallskip \nonumber
&=&\prod_{j=1}^{i}(n^2-(2j-1)^2)
=(-1)^{i/2} \left. a_i \right|_{b_1=n/2}
\end{eqnarray*}
and if $i > n$, then
\begin{eqnarray*}
\prod_{j \in [n;m]_\odd}\left|\frac{i+j}{i-j} \right| \left. a_i\right|_{b_1=0}
(-1)^{\frac{i+n}{2}}
&=&
(-1)^{\frac{i+n}{2}} (i+n-1)!!(i-n-1)!!\\
\smallskip \\
&=&(-1)^{\frac{i+n}{2}+\frac{i-n}{2}} (n-i+1)(n-i+3) \cdots (n-i-1)\\
\smallskip \nonumber
&=&\prod_{j=1}^{i}(n^2-(2j-1)^2)
=(-1)^{i/2} \left. a_i \right|_{b_1=n/2}.
\end{eqnarray*}
Finally, let $n$ be odd and $i$ be even, then
\begin{eqnarray*}
\prod_{j \in [n;m]_\odd}\left|\frac{i+j}{i-j} \right| \left.
a_i\right|_{b_1=0} &=& \frac{(i+n+2m)!!}{(n+2m-i)!!}\\
&=&(n+2m+2-i)(n+2m+4-i) \cdots (n+2m+i)\\
\smallskip \nonumber
&=&\prod_{j=1}^{i}((n+2m+1)^2-(2j-1)^2) =(-1)^{i/2} \left. a_i
\right|_{b_1=(n+2m+1)/2}.
\end{eqnarray*}
\begin{flushright}
$\square$
\end{flushright}

From Lemma {\ref{con.ume.}} we can deduce the following
\begin{Lemma}{\label{con.aux.ham.}}
\begin{eqnarray*}
h_{n,m}(0,b_2)
=\left\{
\begin{array}{lc}
\displaystyle
h_{0,m+\frac{n}{2}}\left(\frac{n}{2},b_2\right),& \mbox{if $n$ is even,}\\
\smallskip \\
\displaystyle
h_{0,\frac{n-1}{2}}\left(m+\frac{n+1}{2},b_2\right),
&\mbox{if $n$ is odd.}
\end{array}
\right.
\end{eqnarray*}
\end{Lemma}
It follows from Lemma {\ref{con.ume.}}, (4.5) and Theorem 1,
that Umemura's polynomials
$U_{m}(b_1,b_2)$ satisfy
a new recurrence relation  with
respect to the first argument $b_1$.
\begin{Theorem}
\begin{eqnarray}
\lefteqn{
U_m(b_1-1,b_2)U_m(b_1+1,b_2)(b_1^2-b_2^2)
}\\ \smallskip \nonumber\\
&=&
(b_1^2-b_2^2)U_m^2(b_1,b_2)+2z^2D_x^2 U_m(b_1,b_2) \circ U_m(b_1,b_2).\nonumber
\end{eqnarray}
\end{Theorem}
Recall that $D_x^2$ denotes the second Hirota derivative.

{\bf Proof.} It follows from Lemma {\ref{con.ume.}} (4.5)
and Theorem 1 that we have
\begin{eqnarray}
\lefteqn{
U_{0,\frac{n-1}{2}}(m+1+\frac{n+1}{2},b_2)U_{0,\frac{n-1}{2}}
(m-1+\frac{n+1}{2},b_2) w^2 \bar{b}_{n+2m+2}
}\\ \smallskip \nonumber\\
&=&
\bar{b}_{n+2m+2} w^2 \left(U_{0,\frac{n-1}{2}}(m+\frac{n+1}{2},b_2)\right)^2
\nonumber \\
\smallskip \nonumber \\
&&+8z^2 w^2 D_x^2 \  U_{0, \frac{n-1}{2}}(m+\frac{n+1}{2},b_2) \circ
U_{0, \frac{n-1}{2}}(m+\frac{n+1}{2},b_2). \nonumber
\end{eqnarray}
We regard the recurrence relation of Theorem 2 as an algebraic
equation with respect to the variable $b_1$. By the above identity
(4.8), this algebraic equation has infinitely many solutions. The
proof of Theorem 2 is finished.

\hfill$\square$
\begin{Corollary}
\begin{eqnarray}
U_{m-1}(b_1,b_2) U_{m+1}(b_1,b_2)
&-& 4 w^2 (b_1^2-b_2^2) U_m(b_1-1,b_2) U_m(b_1+1,b_2) \\
\smallskip \nonumber\\
&=&
\left(
\bar{a}_{2m+2}
\right)
U_m^2(b_1,b_2).\nonumber \\
\smallskip \nonumber \\
U_{m-1}(b_1,b_2) U_{m+1}(b_1,b_2)
&-& 4 z^2 (b_1^2-b_2^2) U_m(b_1,b_2-1) U_m(b_1,b_2+1) \\
\smallskip \nonumber\\
&=& \left( \bar{b}_{2m+2} \right) U_m^2(b_1,b_2), \nonumber
\\\nonumber
\end{eqnarray}
\begin{eqnarray}
U_{m-1}(b_1,b_2) U_{m+1}(b_1,b_2)
&-&   \bar{b}_{2m+2} w^2 U_m(b_1-1,b_2) U_m(b_1+1,b_2)\\
\smallskip \nonumber\\
&+&   \bar{a}_{2m+2} z^2 U_m(b_1,b_2-1) U_m(b_1,b_2+1)
=0.\nonumber
\end{eqnarray}
\end{Corollary}
Let us define functions $X(k,l,m)$ by the following recurrence
relations:
\begin{eqnarray*}
X(k,l,m-1)X(k,l,m+1)&=&X(k-1,l,m)X(k+1,l,m)(4(b_2+l)^2-(2m+1)^2)w^2\\
\smallskip \\
&=&X(k,l-1,m)X(k,l+1,m)(-4(b_1+k)^2+(2m+1)^2)z^2,\\
\end{eqnarray*}
with initial conditions $X(0,0,m)=X(0,1,m)=X(1,0,m)=X(1,1,m)=1$.
To solve these recurrence relations, let us introduce the
following functions.
\begin{eqnarray*}
&&\bar{Y}_{(l,m)}=(4(b_2+l)^2-(2m+1)^2)w^2, \ \
\bar{Z}_{(k,m)}=(-4(b_1+k)^2+(2m+1)^2)z^2, \cr
&&Y_{(l,m)}(n)=\prod_{j=1}^{n}\bar{Y}_{(l,m-n-1+2j)}, ~~~~~~~~~ \
\
\ Z_{(k,m)}(n)=\prod_{j=1}^{n}\bar{Z}_{(k,m-n-1+2j)}.
\end{eqnarray*}
With these notation the explicit formula for $X(k,l,m)$ looks as follows
\begin{eqnarray*}
X(k,l,m)=1\left/ \left(\prod_{j=1}^{k-1}Y_{(l,m)}(j)
\prod_{j=1}^{l-1} Z_{(k,m)}(j)\right). \right.
\end{eqnarray*}
Finally let us introduce the function $T_{(k,l,m)}$ to be
\begin{eqnarray}
T_{(k,l,m)}:=T_{(k,l,m)}(b_1,b_2;z,w)
=U_m(b_1+k,b_2+l) X(k,l,m)
\end{eqnarray}
\begin{Proposition}
Functions $T_{(k,l,m)}$ satisfy the Hirota--Miwa equation
\begin{eqnarray}
T_{(k-1,l,m)}T_{(k+1,l,m)}
+T_{(k,l-1,m)}T_{(k,l+1,m)}
+T_{(k,l,m-1)}T_{(k,l,m+1)}=0.
\end{eqnarray}
\end{Proposition}
This is a direct consequence of Corollary 2.

\section{Example}
\setcounter{equation}{0}
Let us define function $q_m:=q_m(t)$ by the following formula
\begin{eqnarray*}
q_m-t&=&4U_m^2\left\{
\left(m+\frac{1}{2}\right)t(t-1) \frac{d}{dt} \log U_{m+1}
-
\left(m+\frac{3}{2}\right)t(t-1)\frac{d}{dt} \log U_{m}\right.\\
\smallskip \nonumber\\
&&
\left.
-\frac{1}{2} b_1 b_2
+\frac{1}{4}\left(b_1^2\frac{z}{w}+b_2^2\frac{w}{z}\right)
\right\} \left/\left(U_{m+1}U_{m-1}-(2m+1)^2U_m^2\right). \right.
\end{eqnarray*}
One can check that the function $q_m$ is a solution to both
equations $P_{VI}(b_1,b_2,m+\frac{1}{2},0)$ and
$P_{VI}(b_1,b_2,0,m+\frac{1}{2})$. It follows from Okamoto's
theory \cite{OI-OIV} that the function
\begin{eqnarray}
\bar{h}_{1,m}=t(t-1)\frac{d}{dt} \log U_{m+1} -
\frac{1}{4}\left(b_1^2\frac{z}{w}+b_2^2\frac{w}{z}\right)
+\left( m+\frac{1}{2}\right) q_m -\frac{1}{2}\left( m+ \frac{1}{2}\right)
\end{eqnarray}
is also a solution to $E_{VI}(b_1,b_2,m+\frac{1}{2},1)$.
Based on the latter expression for the function $\bar{h}_{1,m}$, and
using Proposition 4 (iii),
we come to the following
\begin{Proposition}\label{gen.sol.}
If $b_1=0$, then we have
\begin{eqnarray*}
U_{m+1}U_{m-1}-(2m+1)^2U_m^2=\frac{1}{4b_2^2}U_{2,m-1}^2,
\end{eqnarray*}
where $U_m:=U_m(0,b_2)$ is a special case of Umemura's polynomial,
and
$$U_{2,m-1}=U_{2,m-1}(0,b_2).
$$
\end{Proposition}
Proof follows from Corollary 2, (4.10).

\begin{flushleft}
\vskip 0.5cm
Anatol N. Kirillov \\
\vskip 5mm
Graduate School of Mathematics

Nagoya University\\

Chikusa-ku, Nagoya,464-8602, Japan\\

and\\

Steklov Mathematical Institute \\

Fontanka 27, St Petersburg 191011, Russia\\
\begin{verbatim}
E-mail address: kirillov@math.nagoya-u.ac.jp
\end{verbatim}
\vskip 1cm
Makoto Taneda \\
\vskip 5mm
Hida 4-3-87, Kumamoto, 861-5514, Japan\\
\begin{verbatim}
E-mail address: tane@rc4.so-net.ne.jp
\end{verbatim}
\end{flushleft}

\end{document}